\pdfoutput=1
\documentclass[11pt]{article}
\usepackage[margin=1in]{geometry}
\usepackage{graphicx}
\usepackage[utf8]{inputenc} 
\usepackage[T1]{fontenc}   
\usepackage{amsfonts}      
\usepackage{amssymb,amsmath,color}
\usepackage{url}
\usepackage{csquotes}
\usepackage{natbib}
\usepackage[mathcal]{eucal}
\usepackage{empheq}
\usepackage{comment}
\usepackage{hyperref}
\newtheorem{remark}{Remark}

\DeclareMathAlphabet\mathbfcal{OMS}{cmsy}{b}{n}

\DeclareMathOperator\KL{\mathrm{KL}}
\newcommand{\R}{\mathbb{R}}
\newcommand{\Pb}{\mathbb{P}}
\newcommand{\E}{\mathbb{E}}
\newcommand{\BEAS}{\begin{eqnarray*}}
\newcommand{\EEAS}{\end{eqnarray*}}
\newcommand{\BEA}{\begin{eqnarray}}
\newcommand{\EEA}{\end{eqnarray}}
\newcommand{\BAL}{\begin{align}}
\newcommand{\EAL}{\end{align}}
\newcommand{\BEQ}{\begin{equation}}
\newcommand{\EEQ}{\end{equation}}
\newcommand{\BIT}{\begin{itemize}}
\newcommand{\EIT}{\end{itemize}}
\newcommand{\BNUM}{\begin{enumerate}}
\newcommand{\ENUM}{\end{enumerate}}
\newcommand{\BA}{\begin{array}}
\newcommand{\EA}{\end{array}}

\begin{document}

\title{The LQR-Schrödinger Bridge}

\author{%
  Marc Lambert  \\
  INRIA - Ecole Normale Supérieure - PSL Research university \\
  DGA - French Procurement Agency\\
  \texttt{marc.lambert@inria.fr} 
}

 \date{}

\maketitle

\begin{abstract}
We consider the Schrödinger bridge problem in discrete time, where the pathwise cost is replaced by a sum of quadratic functions, taking the form of a pathwise linear quadratic regulator (LQR) cost. This cost comprises potential terms, which act as attractors, and kinetic terms, which control the diffusion of the process. When the two boundary marginals are Gaussian, we show that the LQR-Schrödinger bridge problem can be solved in closed form. We follow the dynamic programming principle, interpreting the Kantorovich potentials as cost-to-go functions. Under the LQR-Gaussian assumption, these potentials can be propagated exactly through backward and forward passes, leading to a system of dual Riccati equations, well known in estimation and control. This system converges rapidly in practice. We then show that the optimal process is Markovian and compute its transition kernel in closed form as well as the Gaussian marginals. Through numerical experiments, we demonstrate that this approach can be used to construct complex, non-homogeneous Gaussian processes with acceleration and loops, given well-chosen attractive potentials. Moreover, this approach allows extending the Bures transport between Gaussian distributions to more complex geometries with negative curvature.
\end{abstract}

\section{INTRODUCTION} \label{section1}
In his 1931 foundational article \cite{Chetrite_2021}, Schrödinger addressed the problem of determining the most likely process given two observations of its marginal distributions: $p_0$ at the initial time and $p_K$ at a subsequent time. Following a maximum entropy approach, he considered as a prior, that the process is a Brownian motion starting from $p_0$ and then correct this hypothesis based on the observed marginal $p_K$. At that time, probability theory had not yet been formally established, yet Schrödinger successfully modeled this problem using a system of coupled equations derived purely from physical considerations. This system of equations, known as the Schrödinger system, has attracted renewed interest in recent decades \cite{CLeonard2014, Gentil17, Peyre2019, Yongxin21} due to its connection with entropy-regularized optimal transport and the Sinkhorn algorithm \cite{Sinkhorn64}, which can be used to solve the optimal transport problem efficiently \cite{Cuturi2013}. 

The Schrödinger system is based on the dual formulation of optimal transport, where two Lagrangian functions $V_0$ and $V_K$, also called Kantorovich potentials \cite{Kantorovich1958},  encode the constraints associated with the two bound marginals  $p_0$ and $p_K$.  The potential $V_K$ can be interpreted as a final cost in control theory, which can be propagated backward via dynamic programming to compute the so-called cost-to-go $V_k$ at intermediate times. Without an initial constraint, a single backward pass is sufficient; however, if we impose $p_0$, then a system of forward–backward equations must be solved. When the reference measure is a Brownian motion, efficient solvers such as the Sinkhorn algorithm exist. The case where the reference is a linear stochastic process is more involved and has received significant attention \cite{Levy90, Beghi96, Yongxin16, Bakolas16, Bunne2022}. In particular, in \cite{Yongxin16}, the Schrödinger bridge problem is shown to be equivalent to a stochastic control problem with a quadratic cost on the control where the solution reduces to a system of forward-backward Riccati equations. 

We propose an alternative formulation of the above problem: instead of specifying the reference measure through a linear stochastic process, we define it via a joint distribution determined by a pathwise LQR cost \cite{Kalman60}. This cost incorporate both a quadratic penalty on the control and a quadratic penalty that enforces proximity to a reference path.  
This setting can be useful not only to steer the covariance of the state to a final value \cite{Hotz85} but also to guide it along a path or a tube \cite{Okamoto19}. We show that the resulting problem reduces to a system of Riccati equations bearing a strong resemblance to the Kalman equations for estimation and control. These equations differ from previous works: while the forward equations in \cite{Yongxin16,Bakolas16} describe the evolution of the state covariance matrix for the closed-loop system, our forward equations govern the evolution of the forward Kantorovich potential. Moreover, the pathwise cost introduces additional coupling in the system. Our results are derived in the discrete-time case and restricted to the setting where the velocity is directly controlled, so that the control takes the form 
$u_k \propto x_{k+1}-x_k$. As in the least action principle, the trajectory will then be fully determined by the cost function. This setting is also common in the mean field game approach of Schrödinger bridge \cite{Benamou19}. Extending the framework to control through a channel for a pathwise LQR cost is left for future work.

\begin{figure}[thpb]
\centering
\includegraphics[scale=0.245]{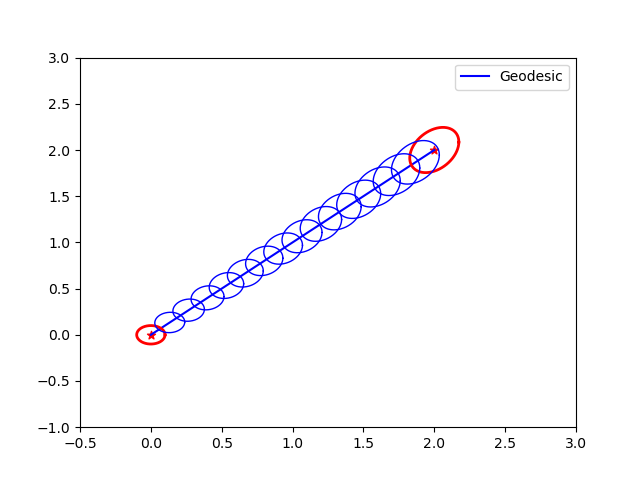}
\includegraphics[scale=0.245]{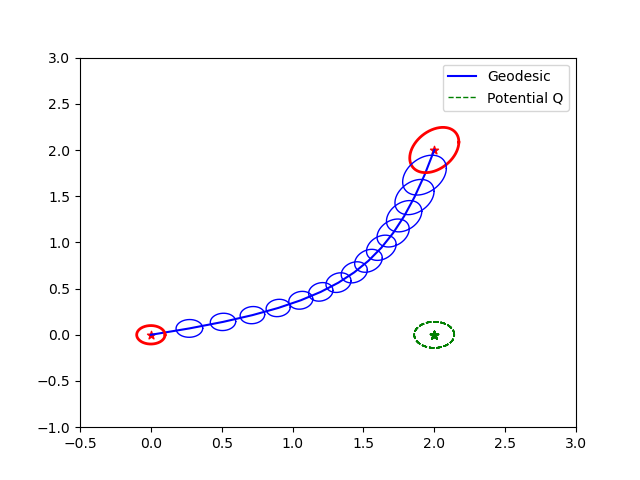}
\includegraphics[scale=0.245]{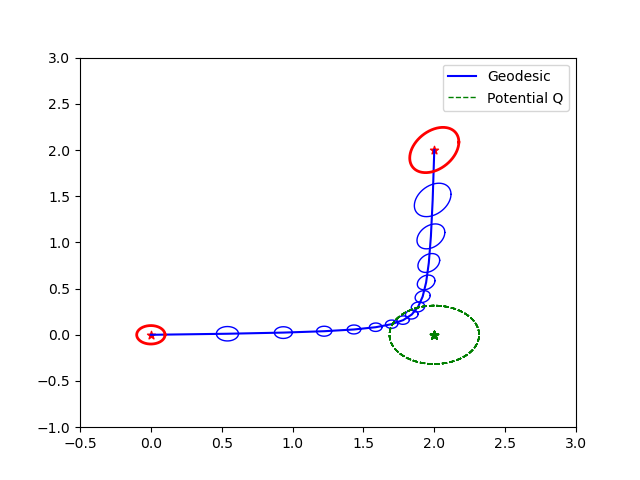}
\includegraphics[scale=0.245]{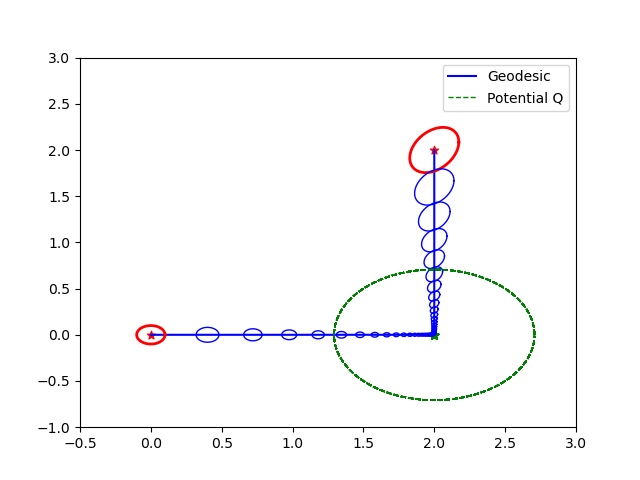}
\caption{As we increase the potential force, shown in green, the Gaussians' Bures geodesics, shown in blue, are deformed and ultimately reduce to a single point where the information is maximal.}
\label{Illustration}
\end{figure}

More precisely, we define our reference measure on the path $x_0,\dots,x_K$ as a Gibbs distribution $\exp(-\frac{1}{\varepsilon} \sum_{k=0}^{K-1} \ell_k(x_k,x_{k+1}))$ where $\varepsilon$ is a temperature parameter and $\ell_k(x_k,x_{k+1})$ is the LQR cost transition which takes the form:
\begin{align*}
&\textbf{The LQR cost}\\
&\boxed{\frac{1}{2}(x_k-x^*_k)^T Q_k (x_k-x^*_k) + \frac{1}{2}(x_{k+1}-x_k)^T R_k (x_{k+1}-x_k).}
\end{align*}
We can reinterpret this cost in the context of stochastic processes as follows: the left-hand term encodes potential energy through the penalization matrix $Q_k$ which attracts   the process around the control points  $x_k^*$. The right-hand term encodes kinetic energy through the penalization matrix  $R_k$, which limits the diffusion of the process. While the kinetic term is consistent with the dynamic formulation of optimal transport \cite{Benamou2000}, the potential term is more unusual and will play an important role in our model.  When considering Gaussian marginals and this quadratic LQR cost, the Kantorovich potential becomes quadratic and can be represented by symmetric positive definite matrices  $P_0^{\oplus},P_K^{\ominus}$ and drift terms $\alpha_0^{\oplus},\alpha_K^{\ominus}$ . We can then show that the forward and backward passes of the Schrödinger system can be computed in closed form, yielding a system of discrete algebraic Riccati equations \cite{Lancaster02} on matrices: 
\begin{align}
&\textbf{The forward-backward discrete Riccati equations} \nonumber\\
&\boxed{P^{\ominus}_{k}=Q_k/\varepsilon+P^{\ominus}_{k+1}- P^{\ominus}_{k+1} ( R_k/\varepsilon+P^{\ominus}_{k+1})^{-1} P^{\ominus}_{k+1}} \label{backRiccati}\\
&\boxed{P^{\oplus^{-1}}_{k+1}=\varepsilon R_k^{-1}+{P^{\oplus}_{k}}^{-1} - {P^{\oplus}_{k}}^{-1} (\varepsilon Q_k^{-1} + {P^{\oplus}_{k}}^{-1} ){P^{\oplus}_{k}}^{-1}.}\label{forwRiccati}
\end{align}
This closely resembles the so-called dual Riccati equations for control (backward) and estimation (forward). 
Using the resulting Kantorovich potentials, we obtain a non-homogeneous Gaussian Markov process governed by the time-varying weighting matrices $Q_k$ and $R_k$. This non-homogeneity offers a versatile modeling framework. In particular, we recover the Bures transport between Gaussians when no potential is active, while activating potentials extends it to more complex geometries problems, as illustrated in Figure \ref{Illustration}.

To compute this solution, we proceed in several steps. First, in Section \ref{Recall-SB} we recall the dual Kantorovich formulation of the Schrödinger problem \cite{CLeonard2001, CLeonard2014} and explain its relation to the optimal transport problem. Next, in Section \ref{Markov-SB}, we consider a structured reference measure with a pairwise loss. In this case, the optimal solution is Markovian, and the non-homogeneous transition kernel can be computed recursively in discrete time \cite{Jamison1975, Beghi96}. We re-derive this result using a simple re-normalization procedure. Our main contribution is given in Section \ref{LQR-SB} where we consider the case in which the pairwise loss function has a LQR quadratic form and assume also  the  marginals are Gaussian. We show the forward and backward passes are given by the Riccati equations on covariance and equations on drifts. In Section \ref{Numerics}, we illustrate the behavior of these optimal Gaussian processes in two dimensions for different types of LQR costs where we activate the potentials $Q_k$ at certain key times $k$ to control the path and shape of the Gaussian process. 


\subsection*{Notations}
We use the convention used in statistical physics and probabilistic graphical models \cite{Wainwright08}. To simplify the exposition,  we assume that all probability measures admit density functions on the $d$-dimensional Euclidean space $\R^d$. This suffices, since we ultimately derive an algorithm for Gaussian distributions, though the results of Sections~\ref{Recall-SB} and \ref{Markov-SB} extend to more general measures, including discrete ones. Moreover, following the graphical model approach, we construct conditional probabilities from joint probabilities in discrete time, i.e., multivariate density functions,  which significantly simplifies the theory of Markov processes. Mathematically speaking, this means we assume that all measures $dp$ considered here are absolutely continuous with respect to the Lebesgue measure
 $dx$ and admit a density function $p(x)$ such that $dp:=p(x)dx$.  $p(x)$ represents the probability distribution $\Pb$ of the random variable $X$: $\Pb(X \in A):=\int_A p(x) dx$ and called without differentiation either a density function or a probability distribution. To simplify notation, we use the same convention as in probabilistic graphical models, and write $p(x)$ and $p(y)$ for distinct probability distributions, using dummy variables to differentiate them. We define a joint probability distribution as $p(x_0,\dots,x_K)$ and assume that the underlying graphical model is a linear chain, such that the index $k$ can be interpreted as time. Given this assumption, we can define the conditional probabilities from this joint probability distribution as follows: $p(x_0,\dots,x_K)=p(x_0) \prod_{k=1}^K p(x_k|x_0,\dots,x_{k-1})$ where the dependencies may simplify if the linear chain is a Markov chain. When $p$ is a normalized probability measure and $f$ not, we use the notation $p(x) \propto f(x) $ to indicate that the equality holds up to a normalization constant such that $p(x)=\frac{1}{Z} f(x)$, where $Z=\int f(x) dx$.
 
 The (negative) entropy is defined by 
   \begin{align*}
   H(dp\|dx):=\int \log \frac{dp}{dx} dp=\int p(x)  \log p(x) dx,
     \end{align*} where $\frac{dp}{dx}$ is the Radon-Nikodym derivative which corresponds to the density function  $p$ with respect to the Lebesgue measure $dx$ here. For brevity, we write $H(dp\|dx):= H(p)$, omitting the reference measure in this case. When considering the relative entropy with respect to a measure $dq=q(x) dx$, which is not simply the Lebesgue measure, we use the Kullback-Leibler (KL) notation:
  \begin{align*}
  &KL(p \| q):=H(dp \| dq)=\int  \log \frac{p(x)}{q(x)} p(x) dx.
  \end{align*}
Finally, given a square-integrable function $V \in \mathrm{L}^2(\R^d)$ and a probability density $p(x)$, we use $\int V(x) p(x) dx := \E_p[V(x)] := \int V dp := \langle V,p \rangle$ interchangeably, with the last one emphasizing the duality product between square-integrable functions and measures. Extending these notations to the joint probability  $p(x_0,x_K)$, we see that a marginal constraint $\int p(x_0,x_K) dx_0=p_K(x_K)$ implies that  $\int V(x_K) (\int p(x_0,x_K) -p_K(x_K) dx_0) dx_K=0$ for any square-integrable functions $V$, or equivalently $\langle V,p-p_K\rangle=0$. We can then see a marginal constraint as a linear constraint in a Hilbert space, which will be useful for constructing the dual Kantorovich problem from a Lagrangian functional. Finally, from the above definitions, we immediately deduce the fundamental relation  $KL(p \| \exp(-V))=H(p)+ \langle V,p \rangle$ which is related to the duality between free energy and entropy. 
 
\section{Recall on Schrödinger bridge} \label{Recall-SB}
\subsection{Static Schrödinger bridge}
Given a cost function $c: \R^d \times \R^d \rightarrow \R,$ the static Schrödinger bridge is given by the  entropic-regularized optimal transport problem from a  measure $p_0$ to a measure $p_K$:
  \begin{align}
\underset{p \in \mathcal{C}(p_0,p_K)}{\min} \E_p[c(x_0,x_K)] + \varepsilon H(p) =\underset{p \in \mathcal{C}(p_0,p_K)}{\min}  \varepsilon \KL \left(p(x_0,x_K) \Big \|  \exp{(-c(x_0,x_K)/\varepsilon)} \right), \label{StaticSB} 
\end{align}
where $\mathcal{C}(p_0,p_K)$ is the set of joint probability distributions (or couplings) on $(x_0,x_K)$ with marginal constraints $(p_0,p_K)$. $\varepsilon$ is a temperature parameter, as $\varepsilon \to 0$, we recover the classical optimal transport problem. 
\subsection{Dynamic Schrödinger bridge}
The above problem can be generalized to the dynamic Schrödinger bridge  if we replace the cost $c(x_0,x_K)$ by a loss  along a path $\ell(x_0,x_1,\dots,x_{K-1},x_K)$ introducing $K-1$ intermediate variables $x_1,\dots,x_{K-1}$:
  \begin{align}
&\underset{p \in \mathcal{P}(p_0,p_K)}{\min} \E_p[\ell(x_0,\dots,x_K)] + \varepsilon H(p)  \label{DynamicSB}\\
=&\underset{p \in \mathcal{P}(p_0,p_K)}{\min} \varepsilon \KL \left(p(x_0,\dots,x_K) \Big \|  \exp{(-\ell(x_0,\dots,x_K)/\varepsilon)} \right) \nonumber\\
=&\underset{p \in \mathcal{P}(p_0,p_K)}{\min} \varepsilon \KL \left(p \Big \|  \exp(-\ell/\varepsilon) \right), \nonumber
\end{align}
where  $\mathcal{P}$ is now the set of joint probabilities on the path  $x_0,\dots,x_K$ with two boundary marginal constraints $(p_0,p_K)$. The right term $\exp(-\ell/\varepsilon)$ is called the reference measure and is not normalized. This lack of normalization plays a fundamental role in the duality between free energy and entropy, as we discuss in the next paragraph. If the reference measure is a Brownian motion, the static and dynamic Schrödinger bridges give the same optimal coupling \cite{Peyre2019}.
 
\subsection{The Kantorovich dual formulation}
We now introduce the Kantorovich dual formulation of the Schrödinger bridge problem, for the derivation in the case of general measures (not just density measures), please refer to  \cite{CLeonard2014,CLeonard2001}. The marginal constraints can be rewritten, for all $V_0,  V_K \in L^2(\R^d)$, as:  $$ \int  V_0 dp= \int  V_0  dp_0; \quad \int   V_K dp = \int  V_K  dp_K,$$ 
where $V_0$ and $V_K$ will play the role of the Lagrange parameters satisfying the linear constraint $\langle V_0, p-p_0 \rangle=0$ and  $\langle V_K ,p-p_K \rangle=0$. These Lagrange parameters are closely related to the so-called dual Kantorovich potentials which play an important role in the theory of optimal transport \cite{Brenier1991, Benamou2000}. Since the KL is convex, there is no duality gap  \cite{CLeonard2001}   and we can switch the min and the max such that the Lagrangian takes the form (where $C_{0,K} :=\int V_0 dp_0+\int V_K dp_K$):
 \begin{align*}
  \underset{p \in \mathcal{P}(p_0,p_K)}{\min} &\KL \left(p \Big \| \exp(-\ell/\varepsilon)  \right)\\
  = \underset{p}{\min} \underset{ V_0,  V_K}{\max} &\KL \left(p \Big \| \exp(-\ell/\varepsilon)  \right) + \langle V_0, p-p_0 \rangle+\langle V_K ,p-p_K \rangle \\
= \underset{ V_0,  V_K}{\max}  \underset{p}{\min} &\KL \left(p \Big \|  \exp(- V_0)  \exp(-\ell/\varepsilon) \exp(-  V_K)  \right) - C_{0,K}.
\end{align*}
The optimal probability is given by $p^* = \frac{1}{Z^*} \exp(- V_0)  \exp(-\ell/\varepsilon) \exp(-  V_K)$ and the optimal value is:
  \begin{align*}
 \underset{ V_0,  V_K}{\max} &- \log Z^* - C_{0,K}\\
 = \underset{ V_0,  V_K}{\max} &- \log \int  \exp(- V_0) \exp(-\ell/\varepsilon)   \exp(-  V_K) - C_{0,K},
\end{align*}
which is the dual free-energy problem \cite{CLeonard2001}.
Finally we rewrite the optimal solution as follow:
\begin{align}
&\boxed{p^*(x_0,\dots,x_K) \propto \phi_0(x_0) r(x_0,\dots,x_K)  \phi_K(x_K),} \label{Solution}
    \end{align}
    where we have introduced the Gibbs potentials  $\phi_0$ and $\phi_K$ supposed normalized and the unnormalized reference measure  $r$ defined as follows: 
\begin{align*}
  &\phi_0(x_0) = \frac{1}{Z_0} \exp(- V_0(x_0)) \\
  &\phi_K(x_K) = \frac{1}{Z_K}  \exp(- V_K(x_K))\\
  &r(x_0,\dots,x_K) =  \exp(-\ell(x_0,\dots,x_K)/\varepsilon).
    \end{align*}

    
    \subsection{The forward-backward  Schrödinger system}
    Reinserting the marginal constraints in the optimal solution \eqref{Solution}, we find that the Gibbs potentials $\phi_0, \phi_K$ must satisfy the Schrödinger system (given here up to a normalization constant):
  \begin{align}
  &\hspace{-0.3cm}  \boxed{p_0(x_0) \propto \phi_0(x_0) \int \phi_K(x_K)  r(x_0,\dots,x_K) dx_1 \dots dx_{K} ,} \label{backward} \\
&\hspace{-0.3cm}  \boxed{p_K(x_K) \propto \phi_K(x_K) \int \phi_0(x_0) r(x_0,\dots,x_K) dx_0 \dots dx_{K-1} .} \label{forward} 
    \end{align} 
When there are no intermediate points and only the two bounds $(x_0,x_K)$, we retrieve the iterative proportional fitting (IPF) algorithm \cite{Csiszar75}. When, in addition, the two marginal measures are discrete, we recover the Sinkhorn algorithm \cite{Sinkhorn64}, which can be computed efficiently \cite{Cuturi2013}. However, when we have a pathwise reference measure with intermediate points, the problem is more involved and requires complex computation of integrals.  In this case, we need to propagate the potential  $\phi_K(x_K)$ backward in equation \eqref{backward}  to update $\phi_0(x_0)$. Indeed, if $\phi^{\ominus}_0(x_0) \propto  \int \phi_K(x_K) r(x_0,\dots,x_K) dx_1 \dots dx_{K}$ is the backward propagated potential then we can compute the initial potential as follows $\phi_0=p_0/\phi^{\ominus}_0$.  The updated potential $\phi_0$ can then be propagated forward in equation \eqref{forward}  to update the terminal potential in the same way $\phi_K=p_K/\phi^{\oplus}_K$. The process is repeated until, hopefully, stabilization of the potentials. Guarantees of convergence may be hard to obtain, but we will see that in the LQR case, fast convergence is observed. This complex forward-backward process simplifies significantly if the reference measure factorizes into a pairwise structure (like Markovian measure). In such cases, the propagation can be computed recursively, allowing us to determine intermediate potentials $\phi_k^{\oplus}$ and  $\phi_k^{\ominus}$. Furthermore, the optimal solution  is necessarily a (non-homogeneous) Markov chain and we can explicitly compute all its transition kernels. This will be the focus of the next section.

\section{The Markovian-Schrödinger bridge} \label{Markov-SB}
\subsection{The pairwise reference measure}
We now refine the setting by considering a loss function that splits into pairwise terms as follows:
\begin{align} 
\ell(x_0,\dots,x_K)=\sum_{k=0}^{K-1} \ell_k(x_k,x_{k+1}). \label{pairwiseLoss}
\end{align} 
\begin{remark}[Connexion with least action principle]
This structured loss can be seen as the discrete counterpart of the action integral $S(x_0)=\int_0^T L(x_t,\dot x_t,t) dt$ where $L$ is the Lagrangian encoding the potential and kinetic energy of a physical system. The action integral $S$ was already proposed as a cost function in optimal transport in the continuous-time setting; leading to the rich theory of  optimal transport on manifolds \cite{villani2008}. From equation \eqref{DynamicSB}, we see also that the Schrödinger bridge problem can be interpreted as a generalization of the principle of least action, with an additive entropic term and bounded constraints. 
\end{remark}

When considering the  loss \eqref{pairwiseLoss}, the unnormalized reference measure factorizes as a pairwise Markov random field \cite{Wainwright08}. Let’s define $r(x_{k+1},x_k)=\exp(-\ell(x_k,x_{k+1})/\varepsilon)$, we have the relation $r(x_0,\dots,x_K)=\exp(-\ell(x_0,\dots,x_K)/\varepsilon)=\prod_{k=0}^{K-1} r(x_{k+1},x_k)$ from which we deduce from equation \eqref{Solution}:
\begin{align*} 
&p^*(x_0,\dots,x_K) \propto \phi_0(x_0) \prod_{k=0}^{K-1} r(x_{k+1},x_k) \phi_K(x_K),
\end{align*} 
and the optimal solution is itself a pairwise Markov random field with bounded constraint.

%
%

From this pairwise structure, we can construct a probability distribution that is a non-homogeneous forward Markov chain $p^*(x_0,\dots,x_K) =p_0(x_0)\prod_{k=0}^{K-1} p^*(x_{k+1}|x_k)$ where the transition kernel $p^*(x_{k+1}|x_k)$ can be constructed from the loss $\ell_k$ as we will show shortly.


\subsection{The optimal Kernel transition}
To compute the optimal transition kernel $p^*(x_{k+1}|x_k)$, we follow the Bellman dynamic programming approach \cite{Bellman2010}. Starting from the final potential $\phi_K(x_K)$, we first compute the last transition $p^*(x_K|x_{K-1})$ and then deduce the form of the next propagated potential $\phi^{\ominus}_{K-1}(x_{K-1})$. We then proceed by recursion to deduce all the subsequent transitions.

Starting from $x_K$, using the pairwise structure of $r$ and the Markov property of $p^*$, we have:
\begin{align*}
p^*(x_K|x_{K-1})&=\phi_0(x_0)  \frac{r(x_0,\dots,x_{K-1})}{p^*(x_0,\dots,x_{K-1})}   r(x_K,x_{K-1})  \phi_K(x_K) \\
&=C r(x_K,x_{K-1})  \phi_K(x_K),
\end{align*}
where we have regrouped all the terms which do not depend on $x_K$ in $C$. Since $p^*(x_K|x_{K-1})$ is a probability distribution, we must have $\int C r(x_K,x_{K-1})  \phi_K(x_K) dx_K=1$ and replacing $C$ by its expression we deduce:
\begin{align*}
p^*(x_0,\dots,x_{K-1})&=\phi_0(x_0) r(x_0,\dots,x_{K-1})   \int r(x_K,x_{K-1})  \phi_K(x_K)  dx_K \\
&= \phi_0(x_0)  r(x_0,\dots,x_{K-1})   \phi^{\ominus}_{K-1}(x_{K-1}),
\end{align*}
where we have introduced the backward-propagated Gibbs potential $$\phi^{\ominus}_{K-1}(x_{K-1}) :=  \int r(x_K,x_{K-1})  \phi_K(x_K)  dx_K=\frac{1}{C}. $$
We can then proceed by recursion to find all the conditionals: 
\begin{align}
&\boxed{p^*(x_{k+1}|x_k) =\frac{1}{\phi^{\ominus}_{k}(x_{k})} r(x_{k+1},x_k) \phi^{\ominus}_{k+1}(x_{k+1})} \label{MarkovSolution}\\ 
&\boxed{ \phi^{\ominus}_k(x_k)= \int r(x_{k+1},x_k)  \phi^{\ominus}_{k+1}(x_{k+1})  dx_{k+1},} \label{RecBackward}
\end{align}
where $r(x_{k+1},x_k)=\exp(-\ell(x_{k+1},x_k))$. 
\begin{remark}[Connexion with Doob’s transform \& RL] Through this re-normalization procedure in discrete time, we recover a Doob’s recursion  \cite{Jamison1975, Beghi96}. Doob's transform is used to condition a reference process to a terminal constraint and is obtained by an importance reweighting of the reference process, where, from equation  \eqref{MarkovSolution}, the importance weight is $\frac{\phi^{\ominus}_{k+1}(x_{k+1})}{\phi^{\ominus}_{k}(x_{k})}$. The recursion \eqref{MarkovSolution}-\eqref{RecBackward} also emerges as a special case of the variational dynamic programming principle introduced in \cite{Lambert24}, specifically when the velocity is controllable, leading to further connections with maximum entropy control \cite{Ziebart2010} and reinforcement learning \cite{Haarnoja2018}. 
\end{remark}

To completely solved the Schrödinger bridge system we need to write also the recursion for the forward pass \eqref{forward}. Using the pairwiswe structure of $r$ it is 
 straightforward, and we obtain:  
\begin{align}
&\boxed{ \phi^{\oplus}_k(x_{k+1})= \int r(x_k,x_{k+1})   \phi^{\oplus}_{k}(x_{k})  dx_{k}.} \label{RecForward}
\end{align}

We have obtained the closed-form expressions for the optimal transition kernel $p^*(x_{k+1}|x_k)$ and the optimal marginals, $p^*(x_k)$, but the computation of the marginalizations in equations \eqref{RecBackward}-\eqref{RecForward} remains intractable in general. We now refine the setting further to find closed-form formulas, considering a quadratic loss and Gaussian marginals. In this case, the Gibbs potentials $\phi_0,\phi_K$ are also Gaussian and preserve their Gaussian structure in both the backward and forward passes. This remarkable property of stability of the quadratic-Gaussian model is at the foundation of the Kalman theory for estimation and control \cite{Kalman60} and is the subject of the next section, where we will consider an LQR loss.

\section{The LQR-Schrödinger bridge}  \label{LQR-SB}
\subsection{Recall on results on Gaussian}
Before delving into the details of the subject, we recall some notations and results on operations with Gaussians in dimension $d$, which will be useful in the sequel. A Gaussian probability distribution with mean $\mu \in \R^d$ and covariance matrix $\Sigma \in \R^{d \times d}$ is denoted by   $q=\mathcal{N}(\mu,\Sigma)$.  The associated Gaussian density function is $q(x)=\mathcal{N}(x|\mu,\Sigma)$. Given these notations, we have the following properties:
\begin{enumerate} 
\item The product of two Gaussian densities $\mathcal{N}(\mu_1,\Sigma_1)$ and $\mathcal{N}(\mu_2,\Sigma_2)$ is, up to a normalization constant $Z$, a Gaussian density, defined by:
\begin{align}
&\mathcal{N}(\mu^*, \Sigma^*)=\frac{1}{Z}  \mathcal{N}(\mu_1,\Sigma_1)\mathcal{N}(\mu_2,\Sigma_2)  \label{ProductGaussians}\\
&\mu^*=\Sigma^*(\Sigma_1^{-1} \mu_1+\Sigma_2^{-1} \mu_2);  \quad \Sigma^*=(\Sigma_1^{-1}+\Sigma_2^{-1})^{-1}; \nonumber \\
& Z =\mathcal{N}(\mu_1|\mu_2,\Sigma_1+\Sigma_2). \label{NormGaussians}
\end{align}
\item The convolution of two Gaussian densities is a Gaussian, defined by:
\begin{align}
&\int \mathcal{N}(x|Az+b,\Sigma_x)\mathcal{N}(z|\mu_z,\Sigma_z) dz \nonumber = \mathcal{N}(A\mu_z+b, \Sigma_x+A\Sigma_zA^T), \label{ConvGaussians}
\end{align}
where  $A \in \R^{d \times d}$ is a matrix and $b \in \R^d$ is a vector.
\end{enumerate}
We will use these results to propagate the Gaussian potentials throught the LQR quadratic cost.

\subsection{The LQR cost with Gaussian marginals}
Inspired by LQR control, we consider a pairwise quadratic loss that mimics the Lagrangian, with a potential energy term and a kinetic energy term:
\begin{align*}
\ell_k(x_k, x_{k+1})&:=\frac{1}{2} (x_k-x^*_k)^T Q_k (x_k-x^*_k) + \frac{1}{2} (x_{k+1}-x_k)^T R_k (x_{k+1}-x_k).
\end{align*}
The potential energy depends on a weighting symmetric positive definite matrix  $Q_k \in d \times d$, which operate as a path-following cost that enforces the process to follow a specific trajectory $x^*_0,\dots,x^*_K$. The kinetic energy depends on a weighting symmetric positive definite matrix $R_k \in d \times d$, which penalizes the diffusion of the process and prevents it from making too big jumps. 

The pairwise refrence measure is then:
\begin{align*}
r(x_k, x_{k+1})&=\exp(-\ell_k(x_k, x_{k+1})/\varepsilon) \propto  \mathcal{N}(x_k|x_k^*,\varepsilon Q_k^{-1})  \mathcal{N}(x_{k+1}|x_k^*,\varepsilon R_k^{-1}).
\end{align*}

Moreover, we suppose that the two bounds marginals are Gaussians distributions given by:
\begin{align*}
& p_0 = \mathcal{N}(\mu_0,\Sigma_0); \quad p_K = \mathcal{N}(\mu_K,\Sigma_K).
\end{align*}

Since $p_0,p_K$ are Gaussians, it is sufficient to choose a Gaussian model for the Gibbs potentials $\phi_0,\phi_K$ to respect the marginals constraints equations \eqref{forward} and \eqref{backward}. 

We thus consider the following Gaussian structure for the unknown Gibbs potentials: 
\begin{align*}
& \phi_0 = \mathcal{N}(\alpha_0,P_0^{-1}); \quad \phi_K = \mathcal{N}(\alpha_K,P_K^{-1}),
\end{align*}
associated to the Kantorovich potentials:
\begin{align*}
& V_0(x_0)=(x_0-\alpha_0)^TP_0(x_0-\alpha_0)\\
& V_K(x_K)=(x_K-\alpha_K)^TP_K(x_K-\alpha_K).
\end{align*}
Finally, at subsequent times $k$, we  note  the forward propagated potential $\phi^{\oplus}_k= \mathcal{N}(\alpha^{\oplus}_k,P^{\oplus^{-1}}_k)$, and the backward propagated potential $\phi^{\ominus}_k= \mathcal{N}(\alpha^{\ominus}_k,P^{\ominus^{-1}}_k)$.

\subsection{The forward-backward Riccati solutions}

\subsubsection{backward potentials }
The recursive update for the backward propagated potentials \eqref{RecBackward} simplifies as:
\begin{align*}
\phi^{\ominus}_k(x_k) &= \int r(x_k,x_{k+1})  \phi^{\ominus}_{k+1}(x_{k+1})  dx_{k+1}\\
&\propto \mathcal{N}(x_k|x^*_k,\varepsilon Q_k^{-1})  \int \mathcal{N}(x_{k+1}|x_k,\varepsilon R_k^{-1})  \mathcal{N}(x_{k+1}|\alpha^\ominus_{k+1},P^{\ominus^{-1}}_{k+1})dx_{k+1}\\
&\propto  \mathcal{N}(x_k|x^*_k,\varepsilon Q_k^{-1})  \mathcal{N}(x_k|\alpha^\ominus_{k+1},\varepsilon R_k^{-1}+P^{\ominus^{-1}}_{k+1}),
\end{align*}
where we have used equation \eqref{NormGaussians}. Using equation \eqref{ProductGaussians} to find the normalized solution,  we finally obtain:
\begin{align*}
&\phi^{\ominus}_k(x_k)=\mathcal{N}(x_k|\alpha^{\ominus}_{k},P^{\ominus^{-1}}_k),
\end{align*}
where we recover equation \eqref{backRiccati} for the covariance matrix:
\begin{align*}
P^{\ominus}_{k}&=Q_k/\varepsilon+(\varepsilon R_k^{-1}+P^{\ominus^{-1}}_{k+1})^{-1}\\
&=Q_k/\varepsilon+P^{\ominus}_{k+1}- P^{\ominus}_{k+1} ( R_k/\varepsilon+P^{\ominus}_{k+1})^{-1} P^{\ominus}_{k+1},
\end{align*}
which comes from the Woodbury formula \cite{Hager89}. 
This update correspond to the LQR backward Riccati equation \cite{Kalman60} in the special case of an identity transition matrix $A_k=\mathbb{I}$. And for the mean we obtain:
\begin{align*}
\alpha^{\ominus}_k&=P^{\ominus^{-1}}_k(Q_k/\varepsilon x^*_k+(\varepsilon R_k^{-1}+P^{\ominus^{-1}}_{k+1})^{-1}\alpha^{\ominus}_{k+1})\\
&=\alpha^{\ominus}_{k+1} + P^{\ominus^{-1}}_kQ_k/\varepsilon (x^*_k-\alpha^{\ominus}_{k+1}),
\end{align*}
which reduce to $\alpha^{\ominus}_k=\alpha^{\ominus}_{k+1}$ if $Q_k=0$ (no bias term). 
\subsubsection{forward potentials }
The recursive update for the forward propagated potentials \eqref{RecForward} simplifies as:
\begin{align*}
\phi^{\oplus}_{k+1}(x_{k+1})&= \int r(x_k,x_{k+1})   \phi^{\oplus}_{k}(x_{k})  dx_k\\
&\propto \int  \mathcal{N}(x_k|x^*_k,\varepsilon Q_k^{-1}) \mathcal{N}(x_k|\alpha^{\oplus}_k,P^{\oplus^{-1}}_k) \mathcal{N}(x_{k+1}|x_k,\varepsilon R_k^{-1})  dx_k\\
&\propto \int  \mathcal{N}(x_k|\alpha^{\oplus}_{k+1},(Q_k/\varepsilon+P^{\oplus}_k)^{-1})\mathcal{N}(x_{k+1}|x_k,\varepsilon R_k^{-1})  dx_k,
\end{align*}
where we have used equation \eqref{ProductGaussians} such that $\alpha^{\oplus}_{k+1}=(Q_k/\varepsilon+P^{\oplus}_k)^{-1}(Q_k/\varepsilon x^*_k+P^{\oplus}_k \alpha^{\oplus}_k)$.
Using equation \eqref{ConvGaussians},  we finally obtain:
\begin{align*}
\phi^{\oplus}_k(x_{k+1})=\mathcal{N}(x_{k+1}|\alpha^{\oplus}_{k+1},P^{\oplus^{-1}}_{k+1}),
\end{align*}
where we recover equation \eqref{forwRiccati} for the covariance matrix:
\begin{align*}
P^{\oplus^{-1}}_{k+1}&=\varepsilon R_k^{-1}+(Q_k/\varepsilon+P^{\oplus}_k)^{-1}\\
&=\varepsilon R_k^{-1}+{P^{\oplus}_{k}}^{-1} - {P^{\oplus}_{k}}^{-1} (\varepsilon Q_k^{-1} + {P^{\oplus}_{k}}^{-1} ){P^{\oplus}_{k}}^{-1},
\end{align*}
which comes again from the Woodbury formula \cite{Hager89}. 
Note that if $Q_k=0$ we may not use this formula. 
\begin{remark}[Duality of Riccati equations] 
Remarkably, we find again a Riccati equation on $P^{-1}$ but where the roles of $R$ and $Q$ are interchanged. This give a new interpretation on the so called Kalman duality where  the inverse  of LQR costs matrixes $R$ and $Q$ can be interpreted as the diffusion matrix of the model noise and \textquote{observation} noise respectively. Indeed $Q_k^{-1}$ can be interepreted as our uncertainty (in terms of covariance) on the prior information. If  $Q_k^{-1}=0$,   then we can predict perfectly where will be the trajectory at $x^*_k$, otherwise we can only estimate a neighborhood around $x^*_k$.  
\end{remark} 

Finally for the drift term we obtain the recursion:
\begin{align*}
\alpha^{\oplus}_{k+1}&=(Q_k/\varepsilon+P^{\oplus}_k)^{-1}Q_k/\varepsilon x^*_k+(Q_k/\varepsilon+P^{\oplus}_k)^{-1} P^{\oplus}_k \alpha^{\oplus}_k,
\end{align*}
which reduce to $\alpha^{\oplus}_{k+1}=\alpha^{\oplus}_k$ if $Q_k=0$ (no bias term). 


\subsubsection{update of potentials }
At the end of the first backward pass, we update the initial potential $\phi_0$ as follows:
\begin{align*}
\phi_0(x_0) =  \frac{p_0(x_0)}{\phi^{\ominus}_0(x_0)}=\frac{\mathcal{N}(\mu_0,\Sigma_0)}{ \mathcal{N}(\alpha^{\ominus}_0,P^{\ominus^{-1}}_0)}. \end{align*}
If $P_0^{\ominus} \prec \Sigma_0^{-1} $ we obtain a well defined Gaussian given by:
\begin{align*}
&\phi_0(x_0)= \mathcal{N}(\alpha_0^*,{P_0^*}^{-1})\\
&P_0^*=\Sigma_0^{-1}-P^\ominus_0\\
&\mu_0^*=(\Sigma_0^{-1}-P^\ominus_0)^{-1}(\Sigma_0^{-1} \mu_0 - P^\ominus_0 \alpha^\ominus_0).
\end{align*}
We then set $\phi^{\oplus}_0(x_0)=\phi_0(x_0)$ and proceed to a forward pass to propagate the updated potential. At the end, we obtain $\phi^{\oplus}_K(x_K)$
which can be used to update the terminal potential $\phi_K(x_K) =  \frac{p_K(x_K)}{\phi^{\oplus}_K(x_K)}$.
We iterate several backward and forward passes until convergence. After convergence, we compute the optimal transition kernel  using an ultimate backward pass.

\subsubsection{optimal transition kernel }
The recursive update for the optimal transition kernel \eqref{MarkovSolution} simplifies as:
\begin{align*}
p^*(x_{k+1}|x_k) &= \frac{1}{\phi_{k}(x_{k})} r(x_k,x_{k+1})   \phi_{k+1}(x_{k+1}) \\
&\propto \mathcal{N}(x_{k+1}|x_k,\varepsilon R_k^{-1})   \mathcal{N}(x_{k+1}|\alpha_{k+1},{P^{\ominus^{-1}}_{k+1}}) \\
&= \mathcal{N}(x_{k+1}|m_k,(R_k/\varepsilon+P^{\ominus}_{k+1})^{-1}),
\end{align*}
where we have used equation \eqref{ProductGaussians} to find the normalized solution and $m_k=(R_k/\varepsilon+P^{\ominus}_{k+1})^{-1}(R_k/\varepsilon. x_k+P^{\ominus}_{k+1}\alpha_{k+1}).$
Finally we obtain:
\begin{align*}
&\boxed{p^*(x_{k+1}|x_k) = \mathcal{N}(x_{k+1}|x_k+\beta_k+K_k x_k,S_k^{-1}),}\\
&\text{with}  \quad S_k=R_k/\varepsilon+P^{\ominus}_{k+1}, \quad K_k=S_k^{-1}R_k/\varepsilon,\\
&\text{and}  \quad \beta_k=S_k^{-1}P^{\ominus}_{k+1}\alpha_{k+1}.
\end{align*}
$K_k$ correspond to the LQR gain when the  transition matrix is the identity and $\beta_k$ is a drift term.

\subsubsection{optimal transport geodesics}
Once we have the optimal transition kernel it is straigthforward to compute the stochastic trajectories. Starting from $x_0 \sim \mathcal{N}(\mu_0,\Sigma_0)$ we have the following random walk in discrete time:
\begin{align*}
&x_{k+1}=x_k+\beta_k+K_k x_k + \nu_k; \quad \nu_k \sim \mathcal{N}(0,S_k^{-1}).
\end{align*}
Computing the mean and the covariance of $x_{k+1}$ we  easily obtain the Gaussian marginals which are given by the recursion:
\begin{align*}
&\boxed{ \mathcal{N}(\mu_{k+1},\Sigma_{k+1})=\mathcal{N}(\beta_k+(\mathbb{I}+K_k) \mu_k,S_k^{-1}+(\mathbb{I}+K_k)^T\Sigma_k(\mathbb{I}+K_k)).}
\end{align*}
\begin{remark}[Connexion with Bures geodesics] 
When $Q_k=0$,  the dynamic Schrödinger bridge \eqref{DynamicSB} is equivalent to the static Schrödinger bridge \eqref{StaticSB} and if $\varepsilon \rightarrow 0$, we recover the classical optimal transport problem. That means these 
sequence of marginals define a geodesics between $\mathcal{N}(\mu_{0},\Sigma_{0})$ and $\mathcal{N}(\mu_{K},\Sigma_{K})$ which matches the Bures-geodesic in this case. When we increase the temperature $\varepsilon$ we obtain the relaxed Bures geodesic \cite{Bunne2022}. More interesting, when $Q_k \succ 0$, we obtain a richer family of  geodesics living on a manifold curved by the potential $Q_k$.
\end{remark}

\section{Numerical Results} \label{Numerics}
We compute now the solution to the LQR-Schrodinger bridge problem in two dimensions ($d=2$) for several types of LQR costs and temperature parameter $\varepsilon$.

\subsection{Initialization}
We provide below the setting used in our simulation and recall the parameters that need to be initialized:

\begin{itemize}
\item The Gaussian marginals mean and covariance $\mathcal{N}(\mu_0,\Sigma_0)$ and $\mathcal{N}(\mu_K,\Sigma_K)$ are shown as red ellipsoids in the plotted figure. The potentials $\phi_0$ and $\phi_K$ are initialized as $\mathcal{N}(\mu_0,\Sigma_0^{-1})$ and $\mathcal{N}(\mu_K,\Sigma_K^{-1})$.
\item The number of points $K+1$ of the path $x_0, \dots, x_K$ is fixed between $15$ and $100$ depending on whether it helps visualization. Note that we use an unscaled discretization step (i.e., $dt=1$), so the value of the loss function depends on the number of points.
\item For the temperature parameter $\varepsilon > 0$ we consider two regimes: the cold regime with $\varepsilon = 0.001$ and the warm regime with $\varepsilon=1$. The number of passes is adapted based on the temperature regime. It is well known that the Sinkhorn algorithm converges slower at low temperatures, and we observe the same behavior in the case of LQR. At the warm temperature or when $K=100$, we observe convergence in just one pass. At the low temperature and for a small number of points $K=15$, we use up to $5$ passes.
\item The diffusion matrix $R_k$ in the LQR cost is assumed fixed to a constant diagonal value $R_k = r \mathbb{I}$ with $r > 0$.
\item The potential matrix is, in most cases, also constant over time and fixed  to a constant diagonal value $Q_k = q \mathbb{I}$ with $q > 0$. The control points $x^*_k$ are chosen to steer the distribution toward arbitrary locations. We show in green the 2D ellipsoid, which corresponds to the potential matrix $Q_k$ centered at $x^*_k$. When no green ellipsoid appears in the plotted figure, there are no control points and $Q_k = 0$. 
\end{itemize}

\subsection{Optimal transport with the LQR geometry}
\begin{figure}[thpb]
\centering
\includegraphics[scale=0.35]{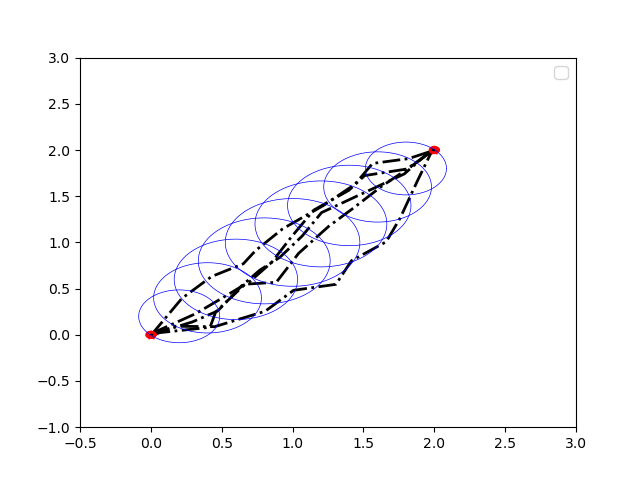}
\includegraphics[scale=0.35]{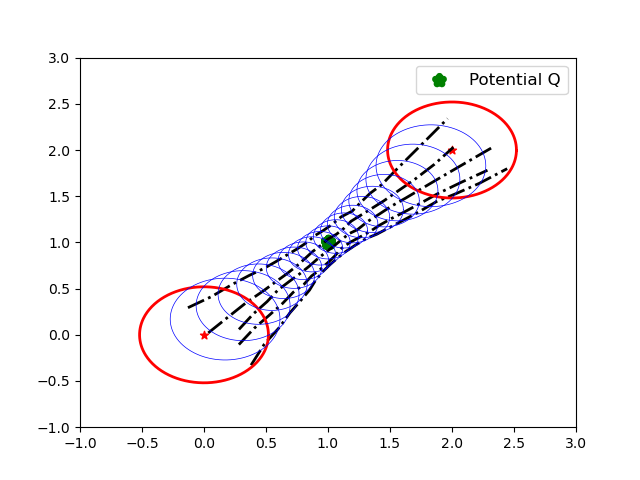}
\caption{The figure on the left shows a standard Brownian bridge between two Dirac marginals, which is obtained when we consider a warm temperature $\varepsilon=1$. There is no potential in the LQR cost, such that $q=0$ and $r=100$. The figure on the right shows an optimal transport between two Gaussians, obtained with a cold temperature $\varepsilon=0.001$. We have added a constant potential in the middle of the path, such that $q=0.3$ and $r=10$.  The covariances of the Gaussian process are shown in blue at $3 \sigma$ here, whereas the sampled trajectories are shown in black. }
\label{Figure1}
\end{figure}

We first illustrate in Figure \ref{Figure1} how the LQR cost can be used to modify the curvature of the optimal transport problem by adjusting the matrices $R$ and $Q$. We consider the optimal transport between two Diracs in the warm regime and between two Gaussians in the cold regime. In the first case, we obtain a Brownian bridge where the paths are stochastic processes (see left plot in Figure \ref{Figure1}). If we decrease the temperature to $0$, these processes shrink into deterministic straight lines,  corresponding to the Euclidean geodesic. In the second case, we replace the Diracs by two Gaussians and consider a low temperature. The optimal transport is given by the Bures transport map, and the samples should form nearly straight lines at this temperature. However, when an attractive potential is added to the path, the resulting geodesics differ from the Bures transport: the geodesics follow a negative curvature  (see right plot in Figure \ref{Figure1}). 

\subsection{Application in stochastic path following and Gaussian steering}
\begin{figure}[thpb]
\centering
\includegraphics[scale=0.245]{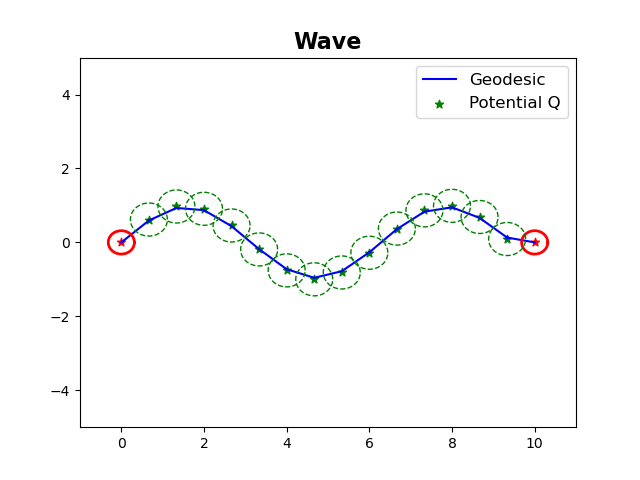}
\includegraphics[scale=0.245]{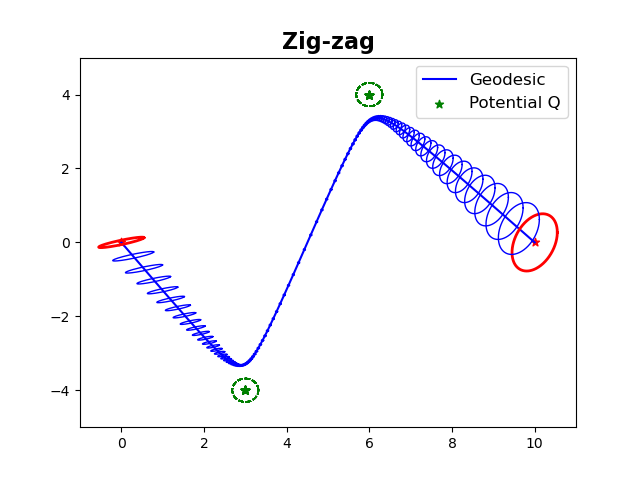}
\includegraphics[scale=0.245]{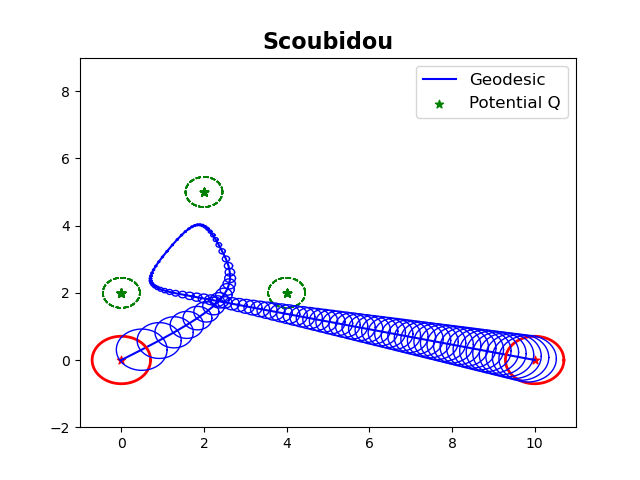}
\includegraphics[scale=0.245]{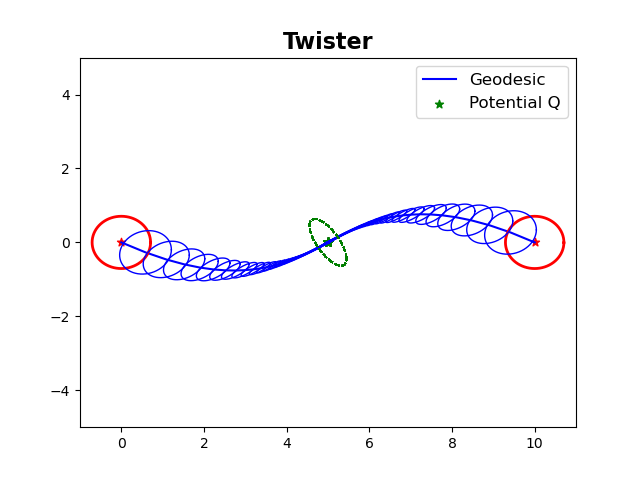}
\caption{Temperature parameter $\varepsilon=0.001$. From left to right:  \textbf{Wave} path following with $r=1, q=10, K=15$, where we also use $K$ waypoints  (green ellipsoid has been rescaled for better visualization in this case); \textbf{Zig-zag} case with $r=10, q=0.1, K=100$, where we use  $2$ waypoints to steer the Gaussian process;\textbf{Scoubidou} case with $r=10, q=0.2, K=100$, where we use  $3$ waypoints; \textbf{Twister} obstacle with $r=10, q=0.2, K=200$, where we use one waypoint with a non-isotropic potential matrix, such that the Gaussian covariances of the marginals are twisted near the potential. The blue ellipsoids show the Gaussian marginals and the green one the potential forces. All the covariances are shown at $1 \sigma$ here. }
\label{Figure2}
\end{figure}

\begin{figure}[thpb]
\centering
\includegraphics[scale=0.245]{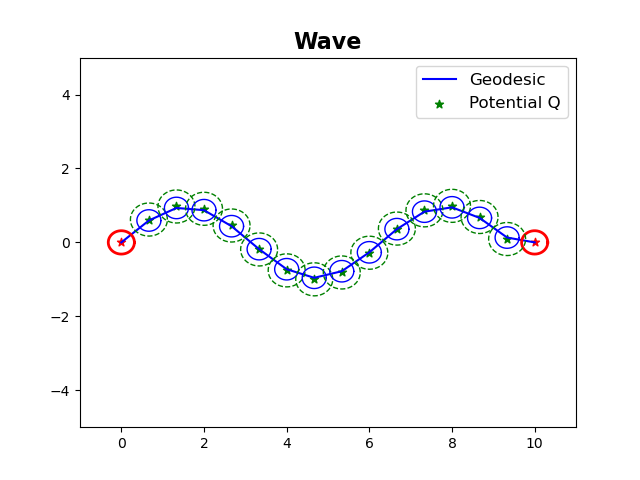}
\includegraphics[scale=0.245]{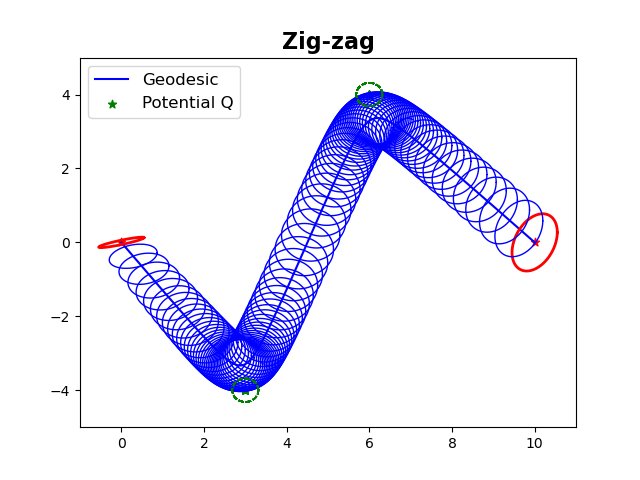}
\includegraphics[scale=0.245]{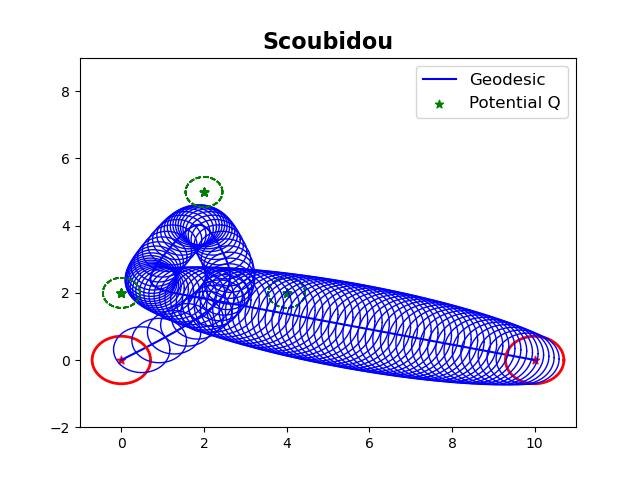}
\includegraphics[scale=0.245]{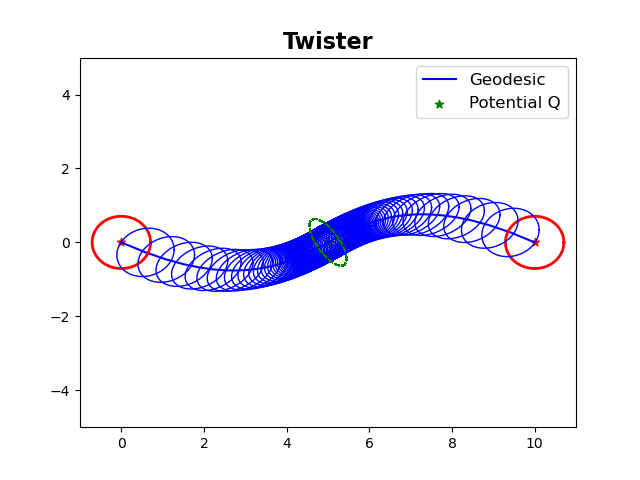}
\caption{Same conditions as in Figure \ref{Figure2}, but with the temperature parameter $\varepsilon=1$. The process is no longer deterministic in the \textbf{Wave} path following and is much more diffusive in the other experiments.}
\label{Figure3}
\end{figure}

Figures \ref{Figure2} and \ref{Figure3} illustrate how the LQR framework can be used to construct complex stochastic processes by modifying the positions of the waypoints or the eigenvalues of the potential covariance matrix. 

This allows the creation of processes that follow a specific path (e.g., the \textquote{Wave} case), avoid obstacles (e.g., the \textquote{Zig-zag} case), form loops (e.g., the \textquote{Scoubidou} case), or accelerate at certain points (e.g., the \textquote{Twister} case). In the cold regime (Figure \ref{Figure2}), the process exhibits little diffusion, and the covariance of the Gaussian decreases due to the potential's influence, in particular in the \textquote{Wave} experiment where the process tends to become deterministic due to the high potential\;value. In the warm regime, the Gaussian process becomes more diffusive and less constrained by the waypoints, as shown in Figure \ref{Figure3}. In all these cases, the potential stabilizes after just one pass.

\textit{The sources of the code are available on Github on the following repository:\\ \url{https://github.com/marc-h-lambert/LQR-SchrodingerBridge}.}

\section{CONCLUSION}
We have integrated ideas from LQR control into the optimal transport framework, leading to a highly effective approach with some surprising implications. Initially, we derive exact closed-form updates for the propagation of the potentials, which take the form of backward and forward Riccati equations that are dual to each other. These equations are inherently discrete-time and can be implemented without approximation, enabling fast convergence of the potentials in practice. Additionally, we compute the optimal transition kernel in closed form, which allows us to generate the associated stochastic process and its marginals. Finally, we observe that a rich geometry emerges when a potential is incorporated into the cost function: by strategically placing attractors, we can control Gaussian distributions, deform them, or reduce their covariance. A natural next step, particularly for control applications, would be to integrate prior information about the process, such as a known passive dynamics, and to restrict the optimization to the control policy $p(u|x)$, linking this work to our earlier research in variational dynamic programming.

\section{ACKNOWLEDGMENTS}
The author gratefully acknowledges C. Léonard for insightful discussions on the dual Kantorovich formulation of the Schrödinger bridge problem, Francis Bach for constructive feedback on the draft that helped emphasize the main ideas, and Amir Taghvaei for the invitation to present this work at the CDC conference. The author also thanks the anonymous reviewer for suggesting relevant references.This work was funded by the French Defence Procurement Agency (DGA) and by the French government under the management of Agence Nationale de la Recherche as part of the “Investissements d’avenir” program, reference ANR-19-P3IA-0001 (PRAIRIE 3IA Institute).
\bibliographystyle{plainnat}
\bibliography{KLcontrol}

\end{document}